\documentclass[12pt]{amsart}
\usepackage[a4paper,bottom=2cm,top=1cm,left=2cm,right=1.3cm]{geometry}
\usepackage[cp1251]{inputenc}
\usepackage[T2A]{fontenc}
\usepackage[english]{babel}
\usepackage{amsfonts,amssymb,amsmath,amsthm,cite}
\usepackage{geometry}
\newtheorem{theorem}{Theorem}

\begin{document}

\title[Multi-interval dissipative Sturm–Liouville problems]{Multi-interval dissipative Sturm–Liouville\\boundary-value problems with distributional coefficients}

\author{Andrii Goriunov}

\address{Institute of Mathematics of National Academy of Sciences of Ukraine, Kyiv, Ukraine}
\email{goriunov@imath.kiev.ua} 

\keywords{Sturm-Liouville operator; multi-interval boundary value problems; distributional coefficients; maximal dissipative extension; completeness of root functions}

\subjclass[2010]{34L40, 34B45}

\begin{abstract}
The paper investigates spectral properties of multi-interval Sturm-Liouville operators with distributional coefficients. 
Constructive descriptions of all self-adjoint and maximal dissipative/accumulative extensions
in terms of boundary conditions are given.  
Sufficient conditions for the resolvents of these operators to be operators of the trace class and for the systems of root functions 
to be complete are found. 
Results of paper are new for one-interval boundary value problems as well.
\end{abstract}

\maketitle 


\textbf{1. Introduction.}

Differential operators, generated by the Sturm-Liouville expression 
\[l(y) = -(py')' + qy, \]
arise in numerous problems of analysis and its applications. 
The classic assumptions on its coefficients are the following: 
\[ q\in C\left([a, b]; \mathbb{R}\right), \qquad 
0< p \in C^1\left([a, b]; \mathbb{R}\right).\]
Principal statements of theory of such operators remain true under more general assumptions
\[ q, 1/p \in L_1\left([a,b], \mathbb{C} \right). \]
However, many problems of mathematical physics require study of differential operators with complex coefficients 
which are Radon measures or even more singular distributions. 
In papers \cite{GM_MFAT_10,Gor-Mikh-ejd-13,MrzSh-16,Teschl-OpMath-13} 
a new approach to investigation of such operators was proposed based on definition of these operators 
as \textit{quasi}-differential, 
which allows also to consider differential operators of higher order  \cite{Gor-Mikh-umzh-11,MrzSh-16}. 

The purpose of this paper is to develop a spectral theory of non-self-adjoint Sturm-Liouville operators,
given on the finite system of bounded intervals under minimal conditions for the regularity of the coefficients.

Multi-interval differential and quasi-differential operators were investigated, 
particularly, in the papers \cite{EvZ_86,EvZ,Sok_06,G_DANU_14}. 
\smallskip

\textbf{2. Preliminary results.}

Let $[a,b]$ --- be a compact interval, $m\in\mathbb{N}$, and let 
$a=a_0<a_1<\dots<a_m=b$ be a partition of  interval $[a,b]$ into $m$ parts. 
Let us consider the space $L_2\left([a,b],\mathbb{C}\right)$  as a direct sum 
$\oplus_{k=1}^mL_2\left([a_{k-1},a_{k}],\mathbb{C}\right)$ which consists of vector functions 
$f= \oplus_{k=1}^m f_k$ such that $f_k \in L_2\left([a_{k-1},a_{k}],\mathbb{C}\right)$. 
 
Let on every interval $(a_{k-1},a_{k})$, $k \in \{1,\dots, m\}$ 
the formal Sturm-Liouville differential expression
\begin{equation}\label{SL_expr_i}
l_k(y)=-(p_k(t)y')'+q_k(t)y+i((r_k(t)y)'+r_k(t)y'),
\end{equation}
be given with coefficients $p_k$, $q_k$ and $r_k$ which satisfy the conditions: 
\begin{equation}\label{SLfull_cond_i}
q_k = Q_k', \quad \frac{1}{\sqrt{|p_k|}}, \frac{Q_k}{\sqrt{|p_k|}}, \frac{r_k}{\sqrt{|p_k|}} \in L_2\left([a_{k-1},a_{k}],\mathbb{C}\right), 
\end{equation}
where the derivatives $Q_k'$ are understood in the sense of distributions. 

Similarly to \cite{MrzSh-16} (see also \cite{GM_MFAT_10,Teschl-OpMath-13}) 
we introduce by the coefficients of the expression \eqref{SL_expr_i} 
on every interval $[a_{k-1},a_{k}]$ the quasi-derivatives in the following way: 
\begin{align*}
D_k^{[0]}y& := y;\\
D_k^{[1]}y& := p_ky' - (Q_k+ir_k)y;\\
D_k^{[2]}y& := (D_k^{[1]}y)'  + \frac{Q_k-ir_k}{p_k}D_k^{[1]}y + \frac{Q_k^2+r_k^2}{p_k}y.
\end{align*} 
Also denote for all $t\in[a_{k-1},a_{k}]$ 
\[\widehat{y}_k(t) = \left( D_k^{[0]}y(t),D_k^{[1]} y(t)\right) \in\mathbb{C}^2.\]
Under assumptions \eqref{SLfull_cond_i} these expressions are Shin-Zettl quasi-derivatives  
(see \cite{EverMarcus-99,Zettl-75}). 
One can easily verify that for the smooth coefficients $p_k$, $q_k$ and $r_k$ the equality  
$l_k(y)~=~-D_k^{[2]}y$ holds. 

Therefore one may correctly define expressions \eqref{SL_expr_i} 
under assumptions \eqref{SLfull_cond_i} as Shin-Zettl quasi-differential expressions: 
$$l_k[y] := - D_k^{[2]}y.$$
The corresponding Shin-Zettl matrices  (see \cite{Zettl-75,EverMarcus-99}) have the form  
\begin{equation}\label{A}
A_k= 
\begin{pmatrix}
\frac{Q+ir}{p}&\frac{1}{p}\\
\,\\
-\frac{Q^2+r^2}{p}&-\frac{Q-ir}{p}
\end{pmatrix}
\in L_1([a,b];\mathbb{C}^{2\times 2}).
\end{equation}

Then on the Hilbert spaces $L_2\left((a_{k-1},a_{k}),\mathbb{C}\right)$ minimal and maximal differential operators 
are defined, which are generated by the quasi-differential expressions $l_k[y]$ (see \cite{EverMarcus-99, Zettl-75}): 
$$ L_{k,1}:y \to
l_k[y],\quad \text{Dom}(L_{k,1}) := 
\left\{y \in L_2 \left| y, D_k^{[1]}y \in AC([a_{k-1},a_k],\mathbb{C}), D_k^{[2]} y \in L_2\right.\right\},$$
 $$ L_{k,0}:y \to
l_k[y],\quad \text{Dom}(L_{k,0}) := 
\left\{y \in \text{Dom}(L_{k,1}) \left| \widehat{y}_k(a_{k-1}) = \widehat{y}_k(a_k) = 0.\right.\right\}.
$$
Results of \cite{Zettl-75,EverMarcus-99} for general Shin-Zettl quasi-differential operators 
together with formula \eqref{A} imply that
operators $L_{k,1}$, $L_{k,0}$ are closed and densely defined 
on the space $L_2\left([a_{k-1},a_{k}],\mathbb{C}\right)$.

In the case where $p_k$, $q_k$ and $r_k$ are real-valued, the operator $L_{k,0}$ is symmetric with
the deficiency index $(2,2)$ and 
\[L_{k,0}^* = L_{k,1},\quad L_{k,1}^* = L_{k,0}.\]
\smallskip

\textbf{3. Dissipative boundary-value problems.}

We consider the space $L_2\left([a,b],\mathbb{C}\right)$ as a direct sum 
$\oplus_{k=1}^mL_2\left([a_{k-1},a_{k}],\mathbb{C}\right)$ which consists of vector functions  
$f= \oplus_{i=1}^m f_i$ such that $f_i \in L_2\left([a_{i-1},a_{i}],\mathbb{C}\right)$.  
In this space $L_2\left([a,b],\mathbb{C}\right)$ we consider maximal and minimal operators 
$L_{\text{max}}=\oplus_{i=1}^mL_{i,1}$ and  
$L_{\text{min}}=\oplus_{i=1}^mL_{i,0}$.

It is easy to see that operators $L_{\text{max}}$, $L_{\text{min}}$ are closed and densely defined on the space 
$L_2\left([a,b],\mathbb{C}\right)$.

Throughout the rest of the paper we assume functions $p_k$, $q_k$ and $r_k$ to be \textit{real-valued} for all $k$ and therefore operators $L_{k,0}$ to be symmetric with the deficiency indices $(2,2)$. 
Then the operator $L_{\text{min}}$ is symmetric with the deficiency index $\left({2m,2m} \right)$ and  
$$L_{\text{min}}^* = L_{\text{max}},\quad L_{\text{max}}^* = L_{\text{min}}.$$
Then naturally arises the problem of describing all its self-adjoint, maximal dissipative and maximal
accumulative extensions in terms of homogeneous boundary conditions of the canonical form. 
For this purpose it is convenient to apply the approach based on the concept of boundary triplets.
It was developed in the papers by Kochubei \cite{Koch-75}, 
see also the monograph \cite{Gorb-book-eng} and the references therein. 

Note that the minimal operator $L_{\text{min}}$ may be not semi-bounded  even 
in the case of a single-interval boundary-value problem 
since the function $p$ may reverse sign.

Recall that a \emph{boundary triplet} of a closed densely defined symmetric operator $T$ 
with equal (finite or infinite) deficiency indices is called a triplet 
$\left(H, \Gamma_1 ,\Gamma_2 \right)$ where $H$ is an auxiliary Hilbert space and 
$\Gamma_1$, $\Gamma_2$ are the linear maps from  $\text{Dom}(T^*)$ into $H$ such that: 
\begin{enumerate}
\item for any $f,g \in \text{Dom}\left(T^*\right)$ there holds 
\[\left(T^*f,g\right)_\mathcal{H} - \left(f,T^*g\right)_\mathcal{H} = \left(\Gamma_1f,\Gamma_2g\right)_H  -
\left(\Gamma_2f,\Gamma_1g\right)_H,\]
\item for any $g_1, g_2 \in H$ there is a vector $f\in \text{Dom}\left(T^*\right)$ 
such that $\Gamma_1 f = g_1$ and $\Gamma_2 f = g_2$.
\end{enumerate}

The definition of the boundary triplet implies that $ f \in \operatorname{Dom} \left( {T} \right)$ 
if and only if $\Gamma_1f = \Gamma_2f = 0$.
A boundary triplet $\left( {H,\Gamma _1 ,\Gamma _2 } \right)$ with $\operatorname{dim} H = n$ 
exists for any symmetric operator $T$ with equal non-zero deficiency indices  $(n, n)$\, $(n \leq \infty)$,
but it is not unique.

\smallskip

For the minimal quasi-differential operators $L_{k,0}$ the boundary triplet is explicitly given 
by the following theorem which follows from the results of \cite{Gor-Mikh-ejd-13}. 
\begin{theorem}\label{Theo_PGZ}
For every $k=1,\dots, m$ the triplet $(\mathbb{C}^{2}, \Gamma_{1,k}, \Gamma_{2,k})$,  
where $\Gamma_{1,k}, \Gamma_{2,k}$ are linear maps  
\[ \Gamma_{1,k}y := \left( D_k^{[1]}y(a_{k-1}+), -D_k^{[1]}y(a_{k}-)\right), \, 
\Gamma_{2,k}y := \left( y(a_{k-1}+), y(a_{k}-)\right), \]
from $\text{Dom}(L_{k,1})$ onto $\mathbb{C}^{2}$ is a boundary triplet for the operator $L_{k,0}$.
\end{theorem}

For the minimal operator $L_{\text{min}}$ in the space $L_2\left([a,b],\mathbb{C}\right)$ 
the boundary triplet is explicitly given by the following theorem. 
\begin{theorem}\label{th_PGZ}
The triplet $(\mathbb{C}^{2m}, \Gamma_{1}, \Gamma_{2})$, 
where $\Gamma_{1}, \Gamma_{2}$ are linear maps 
\begin{equation} \label{PGZ}
 \Gamma_{1}y := \left( \Gamma_{1,1}y, \Gamma_{1,2}y,\dots, \Gamma_{1,m}y\right), \, 
	\Gamma_{2}y := \left( \Gamma_{2,1}y, \Gamma_{2,2}y,\dots, \Gamma_{2,m}y\right), 
\end{equation}
from $\text{Dom}(L_{\text{max}})$ onto $\mathbb{C}^{2m}$ 
is a boundary triplet for the operator $L_{\text{min}}$.
\end{theorem}

Denote by $L_K$ the restriction of operator $L_{\text{max}}$ onto the set of functions 
${y \in \text{Dom}(L_{\text{max}})}$ satisfying the homogeneous boundary condition  
\begin{equation} \label{L_K}
 \left( {K - I} \right)\Gamma_1 y + i\left( {K + I} \right)\Gamma_2 y = 0,
\end{equation}
where $K$ is an arbitrary bounded operator on the space $\mathbb{C}^{2m}$.

Similarly, denote by $L^K$ the restriction of $L_{\text{max}}$ onto the set of functions 
${y \in \text{Dom}(L_{\text{max}})}$ satisfying the homogeneous boundary condition 
\begin{equation} \label{L^K}
 \left( {K - I} \right)\Gamma_1 y - i\left( {K + I} \right)\Gamma_2 y = 0,
\end{equation}
where $K$ is an arbitrary bounded operator on the space $\mathbb{C}^{2m}$.

Theorem \ref{Theo_PGZ} together with \cite[Th.~1.6]{Gorb-book-eng} lead to 
the following description of all self-adjoint, maximal dissipative and maximal accumulative extensions 
of operator $L_{\text{max}}$.

\begin{theorem}\label{th_ext} 
Every $L_K$ with $K$ being a contracting operator in the space $\mathbb{C}^{2m}$,
is a maximal dissipative extension of operator $L_{\text{min}}$.
Similarly every $L^K$ with $K$ being a contracting operator in $\mathbb{C}^{2m}$,
is a maximal accumulative extension of the operator $L_{\text{min}}$.

Conversely, for any maximal dissipative (respectively, maximal accumulative) extension $\widetilde{L}$ 
of the operator $L_{\text{min}}$ there exists the unique contracting operator $K$ such that 
$\widetilde{L} = L_K$\,\, (respectively, $\widetilde{L} = L^K$). 

The extensions $L_K$ and $L^K$ are self-adjoint if and only if $K$ is a unitary operator on $\mathbb{C}^{2m}$. 

The mappings $K\rightarrow L_K$ and $K\rightarrow L^K$ are injective. 
\end{theorem}

All functions from $\text{Dom}(L_{\text{max}})$ together with their first quasi-derivatives 
belong to\\$\oplus_{k=1}^mAC\left([a_{k-1},a_{k}],\mathbb{C}\right)$. 
This implies that following definition is correct.

Denote by $\mathbf{f(t-)}$ the left germ and by $\mathbf{f(t+)}$ the right germ of the continuous function $f$ at point $t$. 
Similarly to the paper \cite{Gor-Mikh-ejd-13} we say that boundary conditions which define the operator 
$L\subset L_{\text{max}}$ 
are called \textit{local}, if for any functions $y\in\text{Dom}(L)$ 
and for any $y_1, \dots, y_m\in \text{Dom}(L_{\text{max}})$ equalities
$\mathbf{y_{j}(a_j-)}=\mathbf{y(a_j-)}$, $\mathbf{y_{j}(a_j+)}=\mathbf{y(a_j+)}$ 
and $\mathbf{y_{j}(a_k-)}= \mathbf{y_{j}(a_k+)}=0$, $k\neq j$  
imply that $y_j\in \text{Dom}(L)$. 
For $j=0$ and $j=m$ we consider only the right and left germs respectively.

The following statement gives a description of extensions $L_K$ and $L^K$ which are given by local boundary conditions. 

\begin{theorem}\label{th_separ}
The boundary conditions \eqref{L_K} and \eqref{L^K} defining extensions $L_K$ and $L^K$ respectively 
are local if and only if the matrix $K$ has the block form
\begin{equation}\label{separ}
K=\begin{pmatrix}K_{a_0}&0&\dots&0\\
0&K_{a_1}&\dots&0\\
0&0&\dots&K_{a_n}
\end{pmatrix},
\end{equation}
where $K_{a_1}$ and $K_{a_n} \in \mathbb{C}$ and other $K_{a_j}\in \mathbb{C}^{2\times 2}$. 
\end{theorem} 

\textbf{4. Generalized resolvents.}

Let us recall that a \textit{generalized resolvent} of a closed symmetric operator $L$ 
in a Hilbert space $\mathcal{H}$ is an operator-valued function $\lambda\mapsto R_\lambda$, 
defined on $\mathbb{C} \setminus \mathbb{R}$ which can be represented as
\[
R_\lambda
f = P^+ \left( L^+ - \lambda I^+\right)^{- 1}f, \quad f \in \mathcal{H},
\]
where $L^+$ is a self-adjoint extension of operator $L$ 
which acts in a certain Hilbert space $\mathcal{H}^+\supset\mathcal{H}$,
$I^+$ is the identity operator on $\mathcal{H}^+$, 
and $P^+$ is the orthogonal projection operator from $\mathcal{H}^+$ onto $\mathcal{H}$.
It is known that an operator-valued function 
$R_\lambda$ 
is a generalized resolvent of a symmetric operator $L$ if and only if 
it can be represented as 
\[\left( R_\lambda f, g \right)_\mathcal{H} = \int_{-\infty}^{+\infty}\frac{d\left(F_\mu f, g\right)}{\mu - \lambda},
\quad f, g \in \mathcal{H},\] 
where $F_\mu$ is a generalized spectral function of the operator $L$.
This implies that the operator-valued function $F_\mu$ has the following properties: 
\begin{enumerate}
\item  For $\mu_2 > \mu_1$ the difference $F_{\mu_2} - F_{\mu_1}$
is a bounded non-negative operator. 
\item $F_{\mu +} = F_\mu$ for any real $\mu$.
\item For any $x \in \mathcal{H}$ following equalities hold: 
\[ \lim\limits_{\mu \rightarrow - \infty}^{}||F_\mu x ||_\mathcal{H} = 0,
\quad \lim\limits_{\mu \rightarrow + \infty}^{} ||{F_\mu x - x} ||_\mathcal{H} = 0.\]
\end{enumerate}

The following theorem provides a parametric description of all generalized resolvents of symmetric operator 
$L_{\text{min}}$ (see also \cite{Brook}). 

\begin{theorem}\label{th_gen_res}
$1)$ 
Every generalized resolvent $R_\lambda$ of the operator $L_{\text{min}}$ in the half-plane $\text{Im}\lambda < 0$
acts by the rule $R_\lambda h = y$, where
$y$ is a solution of the boundary-value problem 
$$ l(y) = \lambda y + h,$$
$$\left( {K(\lambda) - I} \right)\Gamma_{1} f + i\left( {K(\lambda) + I} \right)\Gamma_{2} f = 0.$$
Here $h(x) \in L_2([a,b], \mathbb{C})$ and 
$K(\lambda)$ 
is a $2m\times 2m$ matrix-valued function which is holomorph in the lower half-plane and such that $||K(\lambda)|| \leq 1$.
	
$2)$ In the half-plane $\text{Im}\lambda > 0$ every generalized resolvent of operator $L_{\text{min}}$ 
acts by $R_\lambda h = y$, where $y$ is a solution of the boundary-value problem 
$$ l(y) = \lambda y + h,$$
$$\left( {K(\lambda) - I} \right)\Gamma _{1} f - i\left( {K(\lambda) + I} \right)\Gamma _{2} f = 0.$$
Here $h(x) \in L_2([a,b], \mathbb{C})$ and 
$K(\lambda)$ is a $2m\times 2m$ matrix-valued function which is holomorph
in the upper half-plane and satisfies $||K(\lambda)|| \leq 1$.
\vskip 2mm
This parametrization of the generalized resolvents by the matrix-valued functions $K(\lambda)$ is bijective.
\end{theorem}

\textbf{5. Completeness of system of root vectors.}

Results of the paper \cite{Gor-umzh-15} imply that in the single-interval case 
under the assumptions made and additionally for $r_k=r\equiv 0$ 
the resolvents of the operators $L_K$ and $L^K$ are Hilbert-Schmidt operators. 
This result is amplified and refined by the following theorem.

\begin{theorem}\label{th_root_vectors}
$1)$ The resolvents of the maximal dissipative (maximal accumulative) operators $L_K$ and $L^K$ 
are Hilbert-Schmidt operators.  
	
$2)$ Let $\delta>0$ exist such that for any $k\in \{1,2,\dots,m\}$ 
\[\left\{\frac{1}{p_k}, \frac{Q_k+ir_k}{p_k}\right\} \subset W_2^\delta([a_{k-1},a_{k}],\mathbb{C}). \]
Then the resolvent of the maximal dissipative (maximal accumulative) operator $L_K$ ($L^K$) 
is an operator from the trace class, 
and its system of root functions is complete in the Hilbert space $L_2\left([a,b],\mathbb{C}\right)$. 
\end{theorem}

\end{document}